\documentclass[graybox,envcountsame,envcountsect]{svmult}

\usepackage{mathptmx}       
\usepackage{helvet}         
\usepackage{courier}        
\usepackage{type1cm}        

\usepackage{makeidx}         
\usepackage{graphicx}        
\usepackage{multicol}        
\usepackage[bottom]{footmisc}

\usepackage{url}             
\makeindex             
\usepackage{amsmath}
\usepackage{amssymb}
\usepackage{epsf}
\usepackage{amsfonts}
\usepackage{graphicx}
\usepackage{enumerate}

\newcommand{\BYDEF}{\overset{\mathrm{def}}{=}}

\begin{document}

\title*{Use of approximations of Hamilton-Jacobi-Bellman inequality  for solving periodic optimization problems}
\titlerunning{Approximations of HJB inequality for periodic optimization problems}
\author{  Vladimir Gaitsgory \and
          Ludmila Manic }


\institute{Vladimir Gaitsgory \at Flinders Mathematical Sciences Laboratory, School of Computer Science, Engineering and Mathematics, Flinders University, GPO Box 2100, Adelaide SA 5001, Australia, \email{vladimir.gaitsgory@flinders.edu.au} \and  Ludmila Manic  \at Centre for Industrial and Applicable Mathematics, University of South Australia, Mawson Lakes, SA 5095, Australia \email{Ludmila.Manic@unisa.edu.au} \\ \\ The work was supported
    by the Australian Research Council Discovery-Project Grants DP120100532 and DP130104432}

\maketitle
\abstract{ We show that necessary and sufficient conditions of optimality in periodic optimization  problems can be stated in terms of a solution of the corresponding HJB inequality, the latter being equivalent to a "max-min" type variational problem considered on the space of continuously differentiable functions. We approximate the latter with a maximin problem on a finite dimensional subspace of the space of continuously differentiable functions and show that a solution of this problem (existing under natural controllability conditions) can be used for construction of near optimal controls. We illustrate the construction with a numerical example.}


\section{Introduction and Preliminaries}
\label{sec:1}

\renewcommand{\theequation}{\thechapter.\thesection.\arabic{equation}}
\numberwithin{equation}{section}

%

Consider the control system
\begin{equation}
 y'(t)=f(u(t),y(t)),  \  \   \  \  \ t\in [0,S], \label{1} \end{equation}
where the function $f(u,y):U \times R^m \rightarrow R^m $  is continuous in $(u,y)$ and satisfies Lipschitz conditions in $y$ and where $u(\cdot)$ are controls that are assumed to be Lebesgue measurable and taking values in a given compact metric space $U$.

A pair $(u(\cdot),y(\cdot))$ will be called \emph{admissible} on the interval $[0,S]$ if the equation (\ref{1}) is satisfied for almost all $t\in [0,S]$ and if the following inclusions are valid
\begin{equation} u(t)\in U , \  \  \  \  y(t)\in Y,                 \label{2} \end{equation}
where Y is a given compact subset of $R^m$.

In this paper we will be dealing with an optimal control problem
\begin{equation} \inf_{(u(\cdot),y(\cdot))}\frac{1}{T}\int^T_0 g(u(t),y(t))dt\stackrel{def}{=} G_{per}\ , \label{12} \end{equation}
where $g(u,y):U \times R^m \rightarrow R^1$ is a given continuous function and  \emph{inf} is sought over the length of the time interval $T$ and over the admissible pairs  on $[0,T]$ that satisfy the periodicity condition $y(0)=y(T)$. Problems of this kind
are called {\it periodic optimization} problems. They present a significant mathematical challenge and they have been extensively studied  in the literature
(see, e.g., \cite{Col88}, \cite{FinGaiLeb-08}, \cite{GaiRos06}, \cite{Gra04}).

An important feature of the periodic optimization problems is that they are closely related to so-called long run average optimal control problems. In particular, it can be shown that, under certain conditions, the optimal value of the periodic optimization problem (\ref{12}) is equal to the limit
\begin{equation}   G_{per} = \lim_{S\rightarrow \infty} G(S), \label{15} \end{equation}
where
\begin{equation} G(S)\stackrel{def}{=}\frac{1}{S}\inf_{(u(\cdot),y(\cdot))} \int^S_0 g(u(t),y(t))dt,   \label{3} \end{equation}
 \emph{inf} in (\ref{3}) being sought over all admissible pairs on $[0,S]$ (see, e.g., \cite{Gai91}, \cite{Gai92} and \cite{Gra04}).

Both  problem (\ref{12}) and problem (\ref{3}) can be rewritten in terms of minimization over occupational measures generated by the corresponding admissible pairs. Let $\mathcal{P}(U\times Y)$ stand for the space of probability measures defined on the Borel subsets of $U \times Y$. A measure $\gamma\in \mathcal{P}(U\times Y)$ is called the occupational measure generated by this pair if
it
 satisfies the equation
\begin{equation}  \int_{U\times Y} h(u,y)\gamma(du,dy)=\frac{1}{S} \int^S_0 h(u(t),y(t))dt       \label{4} \end{equation}
for any continuous function $h(u,y):U\times R^m \rightarrow R^1$.

Denote by $\Gamma_{per}$ and $\Gamma (S)$ the sets of occupational measures generated by all periodic admissible pairs and, respectively, by all pairs that are admissible on $[0,S]$. Using these notations, one can equivalently rewrite the problem (\ref{12})   in the form
\begin{equation} \inf_{\gamma \in \Gamma_{per}} \int_{U\times Y}g(u,y) \gamma(du,dy)= G_{per}, \label{5-0} \end{equation}
and the problem (\ref{3}) in the form
\begin{equation} \inf_{\gamma \in \Gamma(S)} \int_{U\times Y}g(u,y) \gamma(du,dy)= G(S). \label{5} \end{equation}
These problems are closely related to the minimization  problem
\begin{equation} \min_{\gamma \in W} \int_{U\times Y} g(u,y)\gamma(du,dy)\stackrel{def}{=}G^*,  \label{6} \end{equation}
where
\begin{equation} W \stackrel{def}{=} \left\{\gamma \in \mathcal{P}(U\times Y): \int_{U\times Y}  \triangledown
\phi(y)^T
f(u,y)\gamma(du,dy)=0, \forall \phi(\cdot) \in C^1\right\}. \label{7} \end{equation}
Note that $C^1$ in the expression above stands for the space of continuously differentiable functions $\phi(\cdot):R^m \rightarrow R^1$, and $\triangledown
\phi(y)
$ is the gradient of $\phi(y)$. Note also that both the objective function in (\ref{6}) and
the constraints in (\ref{7}) are linear in the \lq\lq decision variable" $\gamma$, and, hence the problem (\ref{6}) is one of infinite dimensional (ID) linear programming (see \cite{AndNash87}).

It can be readily shown (see Section 3 in \cite{GaiRos06}) that
\begin{equation} \bar{co}\Gamma_{per}\subset W \ \ \ \ \ \ \Rightarrow \ \ \ \ \ \ \ G^*\leq G_{per},  \label{8-0} \end{equation}
where $ \bar{co}$  in the first expression stands for the closed convex hull. Also, it has been  established  that, under nonrestrictive conditions,
the following relationship is valid  (see Theorem 3.1 in  \cite{FinGaiLeb-08}, Theorem 2.1 in \cite{Gai04} and Proposition 5 in \cite{GaiRos06} as well as related earlier results in  \cite{BhBo96}, \cite{FlemVer-89}, \cite{Her-Her-Tak-96}, \cite{Kur98}, \cite{Vinter} \cite{Rubio85}, \cite{Stock-90})
\begin{equation} \lim_{S\rightarrow \infty} \rho_H (\bar{co}\Gamma(S),W)=0, \label{8} \end{equation}
where $\rho_H (\cdot,\cdot)$ is the Hausdorff metric generated by the weak$^*$ convergence topology (see a precise definition at the end of this section). From (\ref{8})  it, of course, follows that
\begin{equation} \lim_{S\rightarrow \infty}G(S)=G^*  \label{9} \end{equation}
and, if (\ref{15}) is valid, then the latter implies that
\begin{equation}\label{e:14-1}
G_{per} = G^*.  \end{equation}
The validity of the equality (\ref{e:14-1}) is the key assumption  of the present paper. Provided that it is satisfied,
necessary and sufficient optimality conditions for the periodic optimization problem (\ref{12}) are stated in terms
of the Hamilton-Jacobi-Bellman (HJB) inequality, and our main focus is on solving the latter approximately. More specifically (similarly to \cite{GRT}, where infinite horizon optimal control problems with time discounting criterion were considered),
we show that smooth approximate solutions of the HJB inequality exist and that they can be used for the construction of a near optimal control in (\ref{12}).

The paper is organized as follows. In Section 2 that follows this introduction, we define the HJB
inequality and show that it can be used to formulate necessary and sufficient conditions of optimality
for the periodic optimization problem (\ref{12}) (Proposition \ref{Nec-Suf}) . In Section 3, we introduce a variational maximin problem that is equivalent to the HJB inequality  and consider  approximating maximin problems, solutions of which exist (under natural
controllability conditions) and solve the HJB inequality approximately (see Proposition \ref{main-1}).
  In
Section 4, we state a result (Proposition \ref{Main-2})  establishing that solutions of the approximating maxmin problems can be used for construction of a near optimal control in the periodic optimization problem (\ref{12}). In Section 5, we give the proof of this result and in Section 6, we illustrate theoretical developments with a numerical example. The solution of the latter is obtained with a linear programming based algorithm similar to one described
in \cite{GRT}. This algorithm is described in Section 7 and its convergence  (under a non-degeneracy assumption that is less restrictive than one used in \cite{GRT}) is established in Section 8.

 Let us conclude this section with some comments and notations.
 Note,
first of all, that the space ${\cal P} (U\times Y)$ is known to
be compact in the weak$^*$ topology (see, e.g.,
\cite{be_sh} or \cite{Parth}). Being closed, the set $W$
is compact in this topology, and a solution of the problem
(\ref{6})  exists as soon as $W$ is not empty.

Let us  endow the space ${\cal P} (U\times Y)$ with a  metric
$\rho$,
\begin{equation}\label{e:intro-2}
\rho(\gamma', \gamma'') \stackrel{def}{=} \sum^\infty_{j=1} \frac{1}{2^j} \;
\left| \int_{U\times Y} h_j(u,y) \gamma'(du,dy) -  \int_{U\times Y}
h_j(u,y) \gamma''(du,dy)\right| \  ,
\end{equation}
 $\forall \gamma',\gamma''  \in {\cal P} (U\times Y)$, where  $h_j(\cdot), j=1,2,... \ ,$ is a sequence
of Lipschitz continuous functions which is dense in the unit ball
of $C(U\times Y)$ (the space of continuous functions on $U\times
Y$). Note that this metric is consistent with  the weak
convergence topology of ${\cal P} (U\times Y)$. Namely,
 a sequence $\gamma^k \in {\cal P} (U\times Y)$  converges to
$\gamma \in {\cal P} (U\times Y)$ in this metric if and only if
\begin{equation}\label{e:intro-1}
 \lim_{k\rightarrow \infty} \int_{U\times Y} h(u,y) \gamma^k(du,dy) \ = \
 \int_{U\times Y} h(u,y) \gamma(du,dy)
\end{equation}
 for any continuous $h(\cdot)\in C(U\times Y)$.
 Using this metric $\rho$, one can define the $\lq\lq$distance"
$\rho(\gamma , \Gamma)$  between  $ \gamma \in {\cal P} (U\times Y)$ and $\Gamma \subset {\cal P} (U\times Y)$ and the Hausdorff
metric $\rho_H(\Gamma_1, \Gamma_2)$ between $\Gamma_1 \subset
{\cal P} (U\times Y) \ $ and $ \ \Gamma_2 \subset {\cal P} (U\times Y) \ $ as follows:
\begin{equation}\label{e:Housdorff}
\rho(\gamma , \Gamma) \stackrel{def}{=} \inf_{\gamma' \in \Gamma}
\rho(\gamma,\gamma') \ , \ \ \ \ \ \ \ \ \rho_H(\Gamma_1,
\Gamma_2) \stackrel{def}{=} \max \left\{\sup_{\gamma \in \Gamma_1}
\rho(\gamma,\Gamma_2), \sup_{\gamma \in \Gamma_2}
\rho(\gamma,\Gamma_1)\right\} \ .
\end{equation}

\section{Necessary and sufficient conditions of optimality based on the HJB inequality}
\label{sec:2}

The Hamilton-Jacobi-Bellman (HJB) equation for the long run average optimal control problem is written
in the form (see, e.g., Section VII.1.1 in \cite{Bardi97})
\begin{equation} H(\triangledown \psi(y),y)= G^*,   \label{18-0} \end{equation}
where $H(p,y) $ is the Hamiltonian
\begin{equation} H(p,y)\stackrel{def}{=}\min_{u\in U}\{p^T f(u,y)+g(u,y)\}.  \label{17} \end{equation}
The equation (\ref{18-0}) is equivalent to the following two inequalities
\begin{equation} H(\triangledown \psi(y),y)\leq G^*, \ \ \ \   \ \ \  H(\triangledown \psi(y),y)\geq G^*.  \label{18-0-1} \end{equation}
As follows from the result  below, for a characterization of an optimal control in the periodic optimization problem (\ref{12}), it is sufficient to consider
functions
that satisfy only the second inequality in (\ref{18-0-1}), and we will say  that a function $\psi(\cdot)\in C^1$  {\it  is a solution of the HJB inequality on Y}  if
  \begin{equation} H(\triangledown \psi(y),y)\geq G^*,    \  \  \  \forall y \in Y.   \label{18} \end{equation}
Note that the concept of a solution of the  HJB inequality on $Y$ introduced above is essentially the same as that
 of a smooth viscosity subsolution of the HJB equation (\ref{18-0}) considered on the interior of $Y $ (see, e.g., \cite{Bardi97}).

\begin{proposition}\label{Nec-Suf}
Assume that a solution $\psi(\cdot)\in C^1$ of the HJB inequality (\ref{18}) exists. Then  a $T$-periodic admissible  pair $(u(t),y(t))=(u(t+T),y(t+T))$ is optimal in (\ref{12}) and  the equality \begin{equation} G_{per}=G^* \label{1*} \end{equation} is valid if and only if the following relationships are satisfied:
\begin{equation} u(t)=argmin_{u\in U}\{\triangledown \psi(y(t))^T f(u,y(t))+g(u,y(t))\} \  \  a.e. \  t\in [0,T],\label{21} \end{equation}
\begin{equation} H(\triangledown \psi(y(t)),y(t))=G^* \ \ \forall t \in  [0,T]. \label{22} \end{equation}

\end{proposition}

\textbf{Proof}.   Note that from (\ref{17}) and  (\ref{18})  it follows that
\begin{equation} \triangledown \psi(y)^T f(u,y)+ g(u,y)\geq G^*,   \   \   \   \  \forall (u,y)\in U\times Y. \label{22*} \end{equation}
Let us prove the backward implication first. Assume that $G_{per}=G^*$ and  $(u(t),y(t))$ is a solution of the periodic optimization problem (\ref{12}).
That is,
 \begin{equation}\frac{1}{T}\int^T_0 g(u(t),y(t))=G^*. \label{19} \end{equation}
Observe that, due to the periodicity,
\begin{equation} \frac{1}{T}\int^T_0 \triangledown \psi(y(t))^T f(u(t),y(t))dt= \frac{1}{T}\int^T_0 \frac{d(\psi(y(t))}{dt}=\frac{1}{T}(\psi(y(T))-\psi(y(0)))=0.
 \label{28} \end{equation}
From (\ref{19}) and (\ref{28}) it follows that
\begin{equation} \frac{1}{T}\int^T_0 (g(u(t),y(t))+ \triangledown \psi(y(t))^T f(u(t),y(t)))dt=G^*, \label{34} \end{equation}
\begin{equation} \Rightarrow \ \ \ \ \ \frac{1}{T}\int^T_0 (g(u(t),y(t))+ \triangledown \psi(y(t))^T f(u(t),y(t))-G^*)dt=0. \label{34} \end{equation}
By (\ref{22*}), from (\ref{34}) it follows that
\begin{equation} g(u(t),y(t))+\triangledown \psi (y(t))^T f(u(t),y(t))-G^*=0 \  \  a.e. \  t\in [0,T]. \label{30}\end{equation}
 Hence, by (\ref{22*}),
 \begin{equation} (u(t),y(t))\in Argmin_{(u,y)\in U\times Y}\{g(u,y)+\triangledown \psi(y)^T f(u,y)\}\  \  a.e. \  t\in [0,T].  \label{31} \end{equation}
 The latter implies (\ref{21}). Also, by definition of the Hamiltonian (see (\ref{17})), from (\ref{30}) it follows that
 \begin{equation}   H(\triangledown \psi(y(t)),y(t))-G^*\leq 0 \  \  a.e. \  t\in [0,T], \label{31-1}\end{equation}
which (along with the fact that (\ref{18}) is satisfied) prove (\ref{22}).

Let us now prove the forward implication. That is, let us assume that $(u(t),y(t))$ satisfies (\ref{21}) and (\ref{22}), and  show that $(u(t),y(t))$ is an optimal pair and  that $G^*=G_{per}$.
From (\ref{21}) and (\ref{22})  it follows that
\begin{equation} H(\triangledown \psi(y(t)),y(t))= g(u(t),y(t))+\triangledown \psi(y(t))^T f(u(t),y(t))=G^*. \label{32} \end{equation}
By integrating both sides of the above equality  and dividing  by $T$, one obtains
\begin{equation} \frac{1}{T} \int_0^T (g(u(t),y(t))+\triangledown \psi(y(t))^T f(u(t),y(t)))dt=G^*, \label{33} \end{equation}
which, by (\ref{28}),  implies that
\begin{equation} \frac{1}{T} \int_0^T g(u(t),y(t))dt=G^* . \label{33} \end{equation}
 Hence (see (\ref{8-0})), $G_{per}=G^*$ and $(u(t),y(t))$ is optimal.
   $\ \Box$

REMARK. Note that the difference of Proposition \ref{Nec-Suf} from similar results of optimal control theory is  that a solution
of the HJB inequality (rather than that of the HJB equation) is used in the right-hand-side of (\ref{21}), with the relationship (\ref{22}) indicating that the HJB inequality takes the form of the equality on the optimal trajectory.
  Note also that, due to (\ref{18}), the equality (\ref{22}) is equivalent to the inclusion
  \begin{equation} y(t)\in Argmin_{y\in Y}\{H(\triangledown \psi(y),y)\}, \  \  \forall t\in [0,T].  \label{36} \end{equation}

\section{Maximin problem equivalent to the HJB inequality and its approximation}
\label{sec:3}
Consider the following {\it maximin} type problem
  \begin{equation} \sup_{\psi(\cdot) \in C^1} \min_{y\in Y}H(\triangledown \psi(y),y)    \label{19*} \end{equation}
where $sup$ is taken over all continuously differentiable functions.
\begin{proposition}\label{duality}
If the optimal value of the  problem (\ref{19*}) is bounded, then it is equal to the optimal value of the IDLP
problem (\ref{6}). That is,
\begin{equation} \sup_{\psi(\cdot) \in C^1} \min_{y\in Y}H(\triangledown \psi(y),y)=G^*.                \label{201}\end{equation}
\end{proposition}

\textbf{Proof}. As has been shown in \cite{FinGaiLeb-08}, the problem (\ref{19*}) is dual with respect to the IDLP problem
(\ref{6}), and the equality (\ref{201}) follows from the theorem establishing this duality (see Theorem 4.1 in \cite{FinGaiLeb-08}; note that from this theorem it also follows that $supmin$ in (\ref{19*}) is bounded if and only if $W\neq\emptyset$). $\ \Box $

\begin{definition} A function $\psi(\cdot) \in C^1$ will be called a solution of the problem (\ref{19*}) if
\begin{equation} \min_{y\in Y}H(\triangledown \psi(y),y) = G^*.     \label{18-1}\end{equation}
\end{definition}

\begin{proposition}\label{Nec-Suf-Cond}
If $\psi(\cdot) \in C^1$ is a solution of the HJB inequality (\ref{18}), then this $\psi(\cdot)$ is also a solution of the problem (\ref{19*}). Conversely, if $\psi(\cdot)\in C^1$ is a solution of the problem (\ref{19*}), then it also solves the HJB inequality (\ref{18}).
\end{proposition}

\textbf{Proof}. Let  $\psi(\cdot)\in C^1$ be a solution of the HJB inequality
(\ref{18}). By (\ref{19*}) and (\ref{201}), the inequality $\ \min_{y\in Y}H(\triangledown \psi(y),y) > G^*\ $ can not be valid.
Hence, (\ref{18-1}) is true. Conversely, it is obvious that if $\psi(\cdot)$ satisfies (\ref{18-1}), then it satisfies (\ref{18}).
   $\ \Box$

A solution of the maximin problem (\ref{19*})  may not exist, and below we  introduce (following \cite{FinGaiLeb-08}) an \lq\lq approximating" maximin problem. A solution of this problem exists (under non-restrictive conditions) and solves (\ref{19*}) approximately.

Let $\phi_i(\cdot)\in C^1, \  \ i=1,2,...$, be a sequence of functions such that any $\phi(\cdot)\in C^1$ and its gradient are simultaneously approximated by a linear combination of $\phi_i(\cdot), \  \ i=1,2,...,$ and their gradients. An example of such an approximating sequence is the sequence of monomials $y^{i_1}_1 ,..., y^{i_m}_m$,  where  $y_j \ \ (j=1,2,...,m)$ stands for the $j$th component of $y$ and $i_1,..., i_m =0,1,... $ (see e.g.\cite{Lavona86}). Note that it will always be assumed that $\triangledown \phi_i(y), \  \ i=1,2,...,N$ (with $N=1,2,...$), are linearly independent on any open set $Q$. More specifically, it is assumed that, for any $N$, the equality
\begin{equation} \sum^N_{i=1}v_i \triangledown\phi_i(y)=0,  \  \ \forall y \in Q \label{42} \end{equation}
is valid if and only if $v_i=0, \  \ i=1,...,N$.

Define the finite dimensional space $D_N \subset C^1$ by the equation
\begin{equation} D_N \stackrel{def}{=}\left\{\psi(\cdot)\in C^1 : \psi(y)=\sum^N_{i=1}\lambda_i\phi_i(y), \  \ \lambda=(\lambda_i)\in R^N\right\} \label{43} \end{equation}
and consider the maximin problem
\begin{equation} \sup_{\psi(\cdot)\in D_N}\min_{y\in Y}H(\triangledown\psi(y),y)\stackrel{def}{=}\mu^*_N,  \label{44} \end{equation}
which will be referred to as the {\it $N$-approximating maximin problem}.
Note that, due to the definition of the Hamiltonian (\ref{17}), from (\ref{44}) it follows that
\begin{equation} \sup_{\psi(\cdot)\in D_N}\min_{(u,y)\in U\times Y}\{\triangledown\psi(y)^Tf(u,y)+g(u,y)\}=\mu^*_N.  \label{45} \end{equation}

\begin{proposition}\label{N-convergence}
$\mu^*_N$ converges to $G^*$, that is
\begin{equation} \lim_{N\rightarrow \infty} \mu^*_N=G^*.  \label{46} \end{equation}
\end{proposition}

\textbf{Proof}.  It is obvious that, for any $N\geq1$,
\begin{equation} \mu^*_1\leq \mu^*_2\leq...\leq \mu^*_N\leq G^*.                \label{47} \end{equation}
Hence, $\displaystyle \lim_{N\rightarrow\infty}\mu^*_N$ exists, and it is less than or equal to $G^*$. The fact that it is equal to $G^*$ follows from the fact that, for any function $\psi(\cdot)\in C^1$ and for any $\delta > 0$, there exist $N$ large enough and $\psi_\delta (\cdot)\in D_N$ such that
\begin{equation} \max_{y\in Y}\{|\psi(y)-\psi_\delta(y)|+\parallel\triangledown \psi(y)-\triangledown\psi_\delta(y) \parallel\}   \leq \delta .  \label{48}\end{equation}
$\Box$

\begin{definition}\label{def-N-solution}
 A function $\psi(\cdot) \in C^1$ will be called a solution of the $ N$-approximating maximin problem (\ref{44}) if
   \begin{equation}   \min_{y \in Y}\{H(\triangledown \psi(y),y)\}=\mu^*_N  .  \label{50} \end{equation}

\end{definition}
\begin{definition}\label{def-local-controllability}
We shall say that the system (\ref{1}) is {\it locally approximately controllable} on $Y$ if there exists $Y^0 \subset Y$ such that $int (cl Y^0) \neq\emptyset $ (the interior of the closure of $Y^0$ is not empty) and such that any two points in $Y^0$ can be connected by an admissible trajectory.
That is, for any $y_1,y_2\in Y^0$, there exists an admissible pair $(u(t),y(t))$ defined on some interval $[0,S]$ such that $y(0)=y_1$ and $y(S)=y_2$.
\end{definition}

\begin{proposition}\label{main-1} Let the system (\ref{1}) be locally approximately controllable on Y.
   Then, for every $N=1,2,...$, there exists $\lambda^N=(\lambda^N_i)$ such that
   \begin{equation}  \psi^N(y) \stackrel{def}{=}\sum^N_{i=1}\lambda^N_i \phi_i(y)   \label{49} \end{equation}
   is a solution of the  $N$-approximating maximin problem (\ref{44}).
\end{proposition}

\textbf{Proof}. The proof follows from the following two lemmas.

\begin{lemma}\label{important-1}
Assume that, for
\begin{equation} \psi(y)=\sum^N_{i=1} \upsilon_i\phi_i(y),    \label{51}\end{equation}
the inequality
\begin{equation} \triangledown\psi(y)^T f(u,y)\geq 0,   \   \   \forall (u,y)\in U\times Y \label{52} \end{equation}
is valid only if $\upsilon_i=0$,   \  \  $\forall i=1,...,N$.
Then a solution (\ref{49}) of the N-approximating maximin problem (\ref{44}) exists.
\end{lemma}

\begin{lemma}\label{important-2}
If the system (\ref{1}) is locally approximately controllable on $Y$,
then the inequality (\ref{52}) is valid only if $\upsilon_i=0$.
\end{lemma}
\textbf{Proof of Lemma \ref{important-1} }. For any $k=1,2,...$, let $\upsilon^k=(\upsilon^k_i)\in R^N$ be such that the function
\begin{equation} \psi^k(y)\stackrel{def}{=}\sum^N_{i=1} \upsilon^k_i\phi_i(y), \label{53} \end{equation}
satisfies the inequality
\begin{equation} H(\triangledown \psi^k(y),y)\geq \mu^*_N -\frac{1}{k}, \  \  \  \forall y \in Y.                                    \label{54} \end{equation}
Hence,
\begin{equation}   \triangledown \psi^k(y)^Tf(u,y)+g(u,y)\geq \mu^*_N -\frac{1}{k}, \  \  \  \forall(u,y)\in U\times Y.   \label{55} \end{equation}
Let us show that the sequence $\upsilon^k, \ \ k=1,2,...$, is bounded. That is, there exists $\alpha >0$ such that
\begin{equation}   \parallel \upsilon^k \parallel \leq \alpha, \  \  \  k=1,2,...   \ .      \label{56} \end{equation}
Assume that the sequence $\upsilon ^k, \ \ k=1,2,...$, is not bounded. Then there exists a subsequence $\upsilon^{k_l}, \ l=1,2,...$ such that
\begin{equation}   \lim_{l\rightarrow \infty}\parallel \upsilon^{k_l} \parallel =\infty, \    \   \   \lim_{l \rightarrow \infty}\frac{\upsilon^{k_l}}{\parallel\upsilon^{k_l}\parallel} \stackrel{def}{=}\tilde{\upsilon}, \   \   \   \parallel\tilde{\upsilon}\parallel=1.  \label{57} \end{equation}
Dividing (\ref{55}) by $\parallel\upsilon ^{k}\parallel$ and passing to the limit along the subsequence $\{k_l\}$, one can show that
 \begin{equation} \triangledown \tilde{\psi}(y)^T f(u,y)\geq 0,   \   \  \  \forall(u,y)\in U\times Y, \label{58} \end{equation}
 where
\begin{center}{$\tilde{\psi}(y)\stackrel{def}{=}\sum^N_{i=1}\tilde{\upsilon}_i \phi_i(y)$.}\end{center}
Hence, by the assumption of the lemma, $\tilde{\upsilon}=(\tilde{\upsilon_i})=0$, which is in contradiction with (\ref{57}).
Thus, the validity of (\ref{56}) is established.

Due to (\ref{56}), there exists a subsequence $\upsilon^{k_l}, \ l =1,2... ,$ such that there exists a limit
\begin{equation}   \lim_{l \rightarrow \infty} \upsilon^{k_l} \stackrel{def}{=}\upsilon^*.  \label{59} \end{equation} 
Passing to the limit in (\ref{55}) along this subsequence, one obtains
\begin{equation}\triangledown \psi^*(y)^T f(u,y)+ g(u,y)\geq \mu^*_N, \  \  \  \forall(u,y)\in U\times Y, \label{60}\end{equation}
where
\begin{center} $\psi^*(y)\stackrel{def}{=}\sum^N_{i=1}\upsilon^*_i \phi_i(y). $ \end{center}
From (\ref{60}) it follows that
\begin{center} $H(\triangledown \psi^*(y),y)\geq\mu^*_N, \  \  \ \forall y\in Y. $\end{center}
That is, $\psi^*(y)$ is an optimal solution of the N-approximating maximin problem (\ref{44}). $ \ \Box$

\textbf{Proof of Lemma \ref{important-2}}. Assume that
\begin{equation} \psi(y)=\sum^N_{i=1}\upsilon_i\phi_i(y) \label{61} \end{equation}
and the inequality (\ref{52}) is valid. For arbitrary
 $y ^{\prime\ },y^{\prime\prime\ }  \in Y^0$,  there exists an admissible pair  $(u(\cdot), y(\cdot)) $ such that $y(0)=y^{\prime\ }$ and $y(S)=y^{\prime\prime\ }$. From (\ref{52}) it follows that
\begin{center} $\phi(y^{\prime\prime\ })-\phi(y^{\prime\ })=\int^S_0(\triangledown \phi(y(t)))^T f(u(t),y(t))dt \geq 0 \  \  \ \Rightarrow \   \   \  \phi(y^{\prime\prime})\geq \phi(y^{\prime\ }).   $\end{center}
Since $y^{\prime\ }, y^{\prime\prime\ }$ are arbitrary points in $Y_0$, the above inequality allows one to conclude that
$$
\phi(y)=const \ \  \forall y\in Y_0 \ \ \ \ \ \Rightarrow \ \ \ \ \ \phi(y)=const \ \ \forall y\in cl Y^0 ,
$$
the latter implying that $\triangledown\psi(y)=0 \  \  \forall y\in int (clY^0)$ and, consequently leading to the fact that
$\upsilon_i=0, \ i=1,...,N$ (due to the linear independence of $\triangledown\phi_i(y)$). $ \ \Box$

REMARK. Note that from Proposition \ref{N-convergence} it follows that solutions of the $N$-approximating problems (the existence of which is established by Proposition \ref{main-1}))  solve the maximin problem (\ref{19*}) approximately in the sense  that, for any $\delta>0$, there exists $N_\delta$ such that, for any $ N\geq N_\delta$,
\begin{equation}H(\triangledown \psi^N(y),y)\geq G^*-\delta, \  \  \forall y\in Y, \label{63} \end{equation}
where $\psi^N(\cdot)$ is a solution of the $N$-approximating maximin problem (\ref{44}).

\section{Construction of a near optimal control}
\label{sec:4}
In this section, we  assume that a solution $\psi^N(\cdot) $ of the $N$ approximating problem (\ref{44}) exists for all $N$ large enough (see Proposition \ref{main-1}) and we show that, under certain additional assumptions, a control $u^N(y)$ defined as a minimizer of the problem
\begin{equation}\label{u-N} \min_{u\in U} \{\triangledown \psi^N(y)^T f(u,y)+g(u,y)\} \end{equation}
(that is, $u^N(y)=\arg\min_{u\in U} \{\triangledown \psi^N(y)^T f(u,y)+g(u,y)\} $) is near optimal in the periodic optimization problem (\ref{12}). The additional assumptions that we are using to establish this near optimality are as follows.

\textbf{Assumption I}. The equality (\ref{e:14-1})  is valid and  the optimal solution $\gamma^*$ of the IDLP problem (\ref{6}) is unique. Also,  a $T^*$-periodic optimal pair $(u^*(\cdot),y^*(\cdot))$ (that is, the pair that delivers minimum in (\ref{12})) exists.

REMARK. Note that, due to (\ref{e:14-1}), the occupational measure generated by\\ $(u^*(\cdot),y^*(\cdot))$ is an optimal solution of the IDLP problem (\ref{6}). Hence, if $\gamma^*$ is the unique optimal solution of the latter, it will coincide with the occupational measures generated by $(u^*(\cdot),y^*(\cdot))$.

\textbf{Assumption II}. The optimal control $u^*(\cdot):[0,T^*]\rightarrow U$ is piecewise continuous and, at every discontinuity point, $u^*(\cdot)$ is either continuous from the left or it is continuous from the right.

\textbf{Assumption III}.
\\ \textbf{(i)} For almost all $t\in [0,T^*]$, there exists an open ball $Q_t\subset R^m$ centered at $y^*(t)$ such that 
the solution $u^N(y)$ of the problem (\ref{u-N})
is unique for $y\in Q_t $ and that $u^N(\cdot)$ satisfies Lipschitz conditions on $Q_t$ (with a Lipschitz constant being independent of N and t);\\
 \textbf{(ii)} The solution $y^N(\cdot)$ of the system of differential equations
 \begin{equation} y'(t)=f(u^N(y(t)),y(t)), \label{65}  \end{equation}
 which satisfies initial condition $y(0)=y^*(0)$ exists. Moreover this solution is unique and is contained in Y for $t\in [0,T^*]$;\\ \textbf{(iii)} The Lebesgue measure of the set $A_t(N)\stackrel{def}{=} \{ t\in [0,T^*],\ \ \ y^N(t)\notin Q_{t} \}$ tends to zero as $N \rightarrow \infty $. That is,
\begin{equation} \lim_{N\rightarrow\infty}meas\{A_t(N)\}=0. \label{66} \end{equation}

\begin{proposition}\label{Main-2} Let $U$ be a compact subset of $R^n$ and let $f(u,y)$ and $g(u,y)$ be Lipschitz continuous in a neighborhood of $U\times Y$. Also, let the system (\ref{1}) be locally approximately controllable on $Y$ and let Assumptions I, II and III be satisfied. Then
\begin{equation} \lim_{N\rightarrow\infty} u^N(y^N(t))=u^*(t)  \label{67} \end{equation}
for almost all $t\in [0,T^*]$ and
\begin{equation} \max_{t\in [0,T^*]}\parallel y^N(t)-y^*(t)\parallel\leq \nu(N), \ \ \ \ \ \ \ \ \lim_{N\rightarrow\infty}\nu(N)= 0. \label{68} \end{equation}
In addition, if there exists   a $T^*$-periodic solution $\tilde{y}^N(t) $ of the system (\ref{1}) obtained with the control
 $u^N(t)\BYDEF u^N(y^N(t)) $ such that
 \begin{equation}\label{approx-u-per-extra}
\max_{t\in [0,T^*]}\parallel \tilde{y}^N(t)-y^N(t)\parallel\leq \nu_1(N), \ \ \ \ \ \ \ \ \lim_{N\rightarrow\infty}\nu_1(N)= 0,
\end{equation}
 then the pair $(u^N(t),\tilde{y}^N(t))$
 is a near optimal solution of the periodic optimization problem (\ref{12}) in the sense that
\begin{equation} \lim_{N\rightarrow \infty}\frac{1}{T^*} \int^{T^*}_0 g(u^N(t),\tilde{y}^N(t))= G^*. \label{102} \end{equation}
\end{proposition}
The proof of Proposition \ref{Main-2} is given in Section 5.

In conclusion of this section, let us introduce one more assumption, the validity of which implies
the existence of a near optimal periodic admissible pair (see the last part of Proposition \ref{Main-2}).

\textbf{Assumption IV}. The solutions of the system (\ref{1}) obtained with any initial values $y_i, \  i=1,2$ and with any control $u(\cdot)$ satisfy the inequality
\begin{equation}\label{e:stability}\parallel y(t,u(\cdot),y_1)-y(t,u(\cdot),y_2)\parallel\leq \xi(t)\parallel y_1-y_2\parallel,   \ \ \  \ {\rm with} \ \ \ \ \displaystyle \lim_{t\rightarrow \infty}\xi(t)=0. \end{equation}

Note that from Lemma 3.1 in \cite{Gai92} it follows that if Assumption IV is satisfied and if $\xi(T^*)<1$, then the system
\begin{center}$y'(t)=f(u^N(t),y(t))$ \end{center}
(the latter is the system (\ref{1}), in which the control $u^N(t)=u^N(y^N(t)) $ is used)
has a unique $T^*$- periodic solution. Denote this solution as $\tilde{y}^N(T)$.
\begin{proposition}\label{Main-3}
Let Assumptions I,II,III and IV be satisfied and let
\begin{equation}\label{e:xi}
\xi(T^*)<1.
\end{equation}
 Then
\begin{equation}\label{e-new-3} \lim_{N\rightarrow\infty}\max_{t\in[0,T^*]} \parallel \tilde{y}^N(t)-y^N(t) \parallel =0 \end{equation}
and the $T^*$-periodic pair $(u^N(t), \tilde y^N(t)) $ is
a near optimal solution of the periodic optimization problem (\ref{12}) in the sense that (\ref{102}) is valid.
\end{proposition}
The proof is given in Section 5.

\section{Proofs of Propositions \ref{Main-2} and \ref{Main-3}}
\label{sec:5}
Consider the semi-infinite (SI) dimensional LP  problem
\begin{equation} \min_{\gamma\in W_N}\int_{U\times Y}g(u,y)\gamma(du,dy)\stackrel{def}{=}G^*_N, \label{71} \end{equation}
where
 \begin{equation} W_N \stackrel{def}{=} \left\{  \gamma \in P(U\times Y): \ \int_{U\times Y}(\triangledown \phi_i(y)^T f(u,y))\gamma(du,dy)=0, \ \ i=1,...,N \right\} \label{72} \end{equation}
 and $\phi_i(\cdot)$ are as in (\ref{43}).
Note that
 \begin{equation} W_1\supset...\supset W_N\supset W. \label{73} \end{equation}
Consequently, from the fact that $W$ is assumed to be non-empty, it follows that the sets $W_N$, $N=1,2,...$ are not empty. Also (as can be easily seen), the sets $W_N$ are compact in the weak$^*$ topology. Hence, the set of optimal solutions of (\ref{71}) is not empty for any $N=1,2,...$.

\begin{proposition}
The following relationships are valid:
\begin{equation} \lim_{N\rightarrow\infty} \rho_H(W_N,W)=0, \label{74} \end{equation}
\begin{equation} \lim_{N\rightarrow\infty} G^*_N= G^* .      \label{75} \end{equation}
\end{proposition}
\textbf{Proof}. The validity of (\ref{74}) is proved in Proposition 3.5 of \cite{GaiRos06}. The validity of (\ref{75}) follows from (\ref{74}).

\begin{corollary} If the optimal solution $\gamma ^*$ of the problem (\ref{6}) is unique, then for any optimal solution $\gamma^N$ of the problem (\ref{71}) there exists the limit
\begin{equation} \lim_{N\rightarrow\infty} \gamma^N=\gamma^*. \label{76} \end{equation}
\end{corollary}

Note that every extreme point of the optimal solutions set of (\ref{71}) is an extreme point of $W_N$ and that the latter is presented as a convex combination of (no more than $N+1$) 
 Dirac measures (see, e.g., Theorem A.5 in \cite{Rubio85}). That is ,  if $\gamma^N$ is an  extreme point of $W_N$,
which is an optimal solution of (\ref{71}), then there exist
\begin{equation} (u^N_l,y^N_l)\in U\times Y, \  \  \gamma^N_l>0, \  \ l=1,...,K_N\leq N+1; \   \   \sum^{K_N}_{l=1}\gamma^N_l=1 \label{77} \end{equation}
such that
\begin{equation} \gamma^N= \sum^{K_N}_{l=1}\gamma^N_l \delta_{(u^N_l,y^N_l)}, \label{78} \end{equation}
where $\delta_{(u^N_l,y^N_l)}$ is the Dirac measure concentrated at $(u^N_l,y^N_l)$.

The SILP problem (\ref{71}) is related to the $N$-approximating problem (\ref{74}) through the following duality type relationships.
\begin{proposition}
The optimal value of (\ref{71}) and (\ref{44}) are equal
\begin{equation} G^*_N=\mu^*_N. \label{79} \end{equation}
Also, if $\gamma^N$is an optimal solution of (\ref{71}) that allows a representation (\ref{78}) and if $\psi^N(y)=\sum^N_{i=1}\lambda^N_i \phi_i(y)$ is an optimal solution of (\ref{44}), then the concentration points $(u^N_l,y^N_l)$ of the Dirac measures in the expansion (\ref{78}) satisfy the following relationships:
\begin{equation} y^N_l=\arg\min_{y\in Y} \{ H(\triangledown \psi^N(y),y) \}, \label{80} \end{equation}
\begin{equation} u^N_l= \arg\min_{u\in U} \{ \triangledown\psi^N(y^N_l)^Tf(u,y^N_l)+g(u,y^N_l) \}, \  \ l=1,...,K_N. \label{81}
\end{equation}
\end{proposition}

\textbf{Proof}. The validity of (\ref{79}) was proved in Theorem 5.2 (ii) of \cite{FinGaiLeb-07}.
Let us prove (\ref{80}) and (\ref{81}) (note that the argument we are using is similar to that used in \cite{GRT}).

Due to (\ref{79}) and due to the fact that $\psi^N(y)$ is an optimal solution of (\ref{44}) (see (\ref{50})),
\begin{equation} G^*_N=\min_{y\in Y}\{H(\triangledown\psi^N(y),y)\}=\min_{(u,y)\in U\times Y}\{\triangledown\psi^N(y)^T f(u,y)+g(u,y)\}. \label{82} \end{equation}
Also, for any $\gamma\in W_N$,
\begin{center}$\int_{U\times Y}g(u,y)\gamma(du,dy)=\int_{U\times Y} (\triangledown \psi^N(y)^T f(u,y)+ g(u,y))\gamma(du,dy).$\end{center}
Consequently, for $\gamma=\gamma^N$,
\begin{center} $ G^*_N=\int_{U\times Y} g(u,y)\gamma^N(du,dy)=\int_{U\times Y}[g(u,y)+\triangledown\psi^N(y)^T f(u,y)]\gamma^N(du,dy).$\end{center}
Hence, by (\ref{78}),
\begin{equation} G^*_N=\sum^{K_N}_{l=1}\gamma^N_l[g(u^N_l,y^N_l)+\triangledown\psi^N(y^N_l)^T f(u^N_l,y^N_l)]. \label{83} \end{equation}
Since $(u^N_l,y^N_l)\in U\times Y$, from (\ref{82}) and (\ref{83}) it follows that, if $\gamma^N_l>0$, then 
\begin{center}$\displaystyle g(u^N_l,y^N_l)+\triangledown\psi^N (y^N_l)^T f(u^N_l,y^N_l)= \min_{(u,y)\in U\times Y} \{\triangledown \psi^N(y)^T f(u,y)+g(u,y)\}.$\end{center}
That is,
\begin{center}$\displaystyle (u^N_l,y^N_l)=\arg\min_{(u,y)\in U\times Y}\{\triangledown\psi^N(y)^T f(u,y)+ g(u,y)\}.$\end{center}
The latter is equivalent to (\ref{80}) and (\ref{81}).  \  \  \  $\ \Box$

\begin{lemma} Let Assumptions I and II be satisfied and let $\gamma^N$ be an optimal solution of (\ref{71}) that is presented in the form (\ref{78}). Then
\begin{equation} \sup_{t\in[0,T^*]}d((u^*(t),y^*(t)),\Theta_N)\rightarrow 0 \  \  \  \ as\  \  N\rightarrow\infty, \label{84} \end{equation}
where  $\Theta_N\stackrel{def}{=}\{(u^N_l,y^N_l ), \ \ l=1,...,K_N\}$.
\end{lemma}

\textbf{Proof}. Let $\Theta^*\stackrel{def}{=}\{(u,y): (u,y)=(u^*(t),y^*(t))\  {\rm for \ some} \ t\in[0,T^*]\}$, and let $B$ be the open ball in $R^{n+m}$: $B\stackrel{def}{=}\{(u,y):\parallel(u,y)\parallel<1\}$. It is easy to see that Assumption II implies that, for any $(\overline{u},\overline{y})\in cl\Theta^*$ (the closure of $\Theta^*$) and any $r>0$, the set $B_r(\overline{u},\overline{y})\stackrel{def}{=}((\overline{u},\overline{y})+rB)\cap(U\times Y)$ has a nonzero $\gamma^*$-measure. That is,
\begin{equation} \gamma^*(B_r(\overline{u},\overline{y}))>0. \label{85} \end{equation}
In fact, if $(\overline{u},\overline{y})\in cl\Theta^*$, then there exists a sequence $t_i, \ \ i=1,2,...$, such that $(\overline{u},\overline{y})=\lim_{i\rightarrow\infty}(u^*(t_i),y^*(t_i))$, with $(u^*(t_i),y^*(t_i))\in B_r(\overline{u},\overline{y})$ for some $i$ large enough.
Hence, there exists $\alpha >0$ such that $(u^*(t^{\prime\ }),y^*(t^{\prime\ }))\in B_r(\overline{u},\overline{y}),  \  \  \forall t^{\prime\ }\in (t_i-\alpha,t_i]$ if $u^*(\cdot)$ is continuous from the left at $t_i$, and $(u^*(t^{\prime\ }),y^*(t^{\prime\ }))\in B_r(\overline{u},\overline{y}), \  \  \forall t^{\prime\ }\in [t_i,t_i+\alpha)$ if $u^*(\cdot)$ is continuous form the right at $t_i$. Since $\gamma^*$ is the occupational measure generated by the  pair $(u^*(t),y^*(t)) $ (see Remark after Assumption I), the latter   implies (\ref{85}).

Assume now the statement of the lemma is not valid. Then there exist a number $r>0$ and sequences: $N_i, (u^*(t_i),y^*(t_i))\stackrel{def}{=}(u_i,y_i) \in \Theta^*, i=1,2,...$ , with
\begin{center} $\displaystyle \lim_{i\rightarrow \infty}(u_i,y_i)=(\overline{u},\overline{y})\in cl\Theta, \  \  \  \lim_{i\rightarrow\infty} N_i=\infty,$ \end{center}
such that
\begin{equation} d((u_i,y_i),\Theta^{N_i})\geq 2r \  \  \  \Rightarrow \  \  \ d((\overline{u},\overline{y}),\Theta^{N_i})\geq r, \  \ i\geq i_0,\label{88}\end{equation}
where $d((u,y),Q)$ stands for the distance between a point $(u,y)\in U\times Y$ and a set $\displaystyle Q\subset U\times Y: d((u,y),Q)\stackrel{def}{=}\inf_{(u^{\prime },y^{\prime })\in Q}\{ \parallel(u,y)-(u^{\prime\ },y^{\prime\ })\parallel \} $. The second inequality in (\ref{88}) implies that
\begin{center} $(u^{N_i}_l,y^{N_i}_l)\notin B_r(\overline{u},\overline{y}), \  \  \  l=1,...,K_{N_i}, \  \ i\geq i_0.$ \end{center}
By (\ref{77}), the latter implies that
\begin{equation} \gamma^{N_i}(B_r(\overline{u},\overline{y}))=0. \label{89} \end{equation}
From (\ref{76}) it follows that
\begin{center}$\displaystyle \lim_{i\rightarrow\infty} \rho(\gamma^{N_i},\gamma^*)=0.$ \end{center}
Consequently (see, e.g., Theorem 2.1 in \cite{Bill68}),
\begin{center}$\displaystyle 0= \lim_{i\rightarrow\infty} \gamma^{N_i}(B_r(\overline{u},\overline{y})\geq\gamma^*(B_r(\overline{u},\overline{y})))$. \end{center}
The latter contradicts (\ref{85}) and thus proves the lemma.  $\ \Box$

\textbf{Proof of Proposition \ref{Main-2}}. Let $t\in [0,T^*]$ be such that $u^N(\cdot)$ is Lipschitz continuous on $Q_t$. By (\ref{84}), there exists $(u^N_{l_N},y^N_{l_N})\in \Theta_N$ such that
\begin{equation}\label{e-extra103}
\lim_{N\rightarrow\infty}||(u^N_{l_N},y^N_{l_N})- (u^*(t),y^*(t))||=0,
\end{equation}
the latter implying, in particular, that $y^N_{l_N}\in Q_t $ for $N$ large enough.
 Due to (\ref{81}),
\begin{equation}u^N_{l_N}=u^N(y^N_{l_N}). \label{86} \end{equation}
Hence,
\begin{center}$\parallel u^*(t) - u^N(y^*(t)\parallel \leq \parallel u^*(t)-u^N_{l_N}\parallel+\parallel u^N(y^N_{l_N})-u^N(y^*(t))\parallel$\end{center}
\begin{equation} \leq \parallel u^*(t)-u^N_{l_N}\parallel+L\parallel y^*(t)-y^N_{l_N}\parallel, \label{87}\end{equation}
where L is a Lipschitz constant of $u^N(\cdot)$. From (\ref{e-extra103})  it now follows that
\begin{equation} \lim_{N\rightarrow\infty} u^N(y^*(t))=u^*(t). \label{91} \end{equation}
By Assumption III, the same argument is applicable for almost all $t\in [0,T^*]$. This proves the convergence (\ref{91})
for almost all $t\in [0,T^*]$.

Taking an arbitrary $t\in[0,T^*]$ and subtracting the equation
\begin{equation} y^*(t)=y_0+\int^t_0 f(u^*(t^{\prime\ }), y^*(t^{\prime\ }))dt^{\prime\ }  \label{92} \end{equation}
from the equation
\begin{equation} y^N(t)=y_0+\int^t_0 f(u^N(y^N(t^{\prime\ })), y^N(t^{\prime\ }))dt^{\prime\ },  \label{93} \end{equation}
one obtains
\begin{center} $\parallel y^N(t)-y^*(t)\parallel \leq\int^t_0 \parallel f(u^N(y^N(t^{\prime\ })),y^N(t^{\prime\ }))-f(u^*(t^{\prime\ }),y^*(t^{\prime\ }))\parallel dt^{\prime\ }$\end{center}
\begin{center} $\leq\int^t_0\parallel f(u^N(y^N(t^{\prime\ })),y^N(t^{\prime\ }))-f(u^N(y^*(t^{\prime\ })),y^*(t^{\prime\ }))\parallel dt^{\prime\ }$  \end{center}
\begin{equation}  +\int^t_0\parallel f(u^N(y^*(t^{\prime\ })),y^*(t^{\prime\ }))-f(u^*(t^{\prime\ }),y^*(t^{\prime\ }))\parallel dt^{\prime\ }. \label{94} \end{equation}
It is easy to see that
\begin{center} $\int^t_0\parallel f(u^N(y^N(t^{\prime\ })),y^N(t^{\prime\ }))-f(u^N(y^*(t^{\prime\ })),y^*(t^{\prime\ }))\parallel dt^{\prime\ }$ \end{center}
\begin{center} $\leq \int_{t^{\prime\ }\notin A_t(N)} \parallel f(u^N(y^N(t^{\prime\ })),y^N(t^{\prime\ }))-f(u^N(y^*(t^{\prime\ })),y^*(t^{\prime\ }))\parallel dt^{\prime\ }$ \end{center}
\begin{center} $+ \int_{t^{\prime\ }\in A_t(N)}[ \parallel f(u^N(y^N(t^{\prime\ })),y^N(t^{\prime\ }))\parallel+\parallel f(u^N(y^*(t^{\prime\ })),y^*(t^{\prime\ }))\parallel] dt^{\prime\ }$ \end{center}
\begin{equation} \leq L_1 \int^t_0 \parallel y^N(t^{\prime\ })-y^*(t^{\prime\ }) \parallel dt^{\prime\ }+L_2 \ meas\{A_t(N)\}, \label{95} \end{equation}
where $L_1$ is a constant defined (in an obvious way) by Lipschitz constants of $f(\cdot,\cdot)$ and $u^N(\cdot)$, and $\displaystyle L_2\stackrel{def}{=}2 \max_{(u,y)\in U\times Y}\{\parallel f(u,y)\parallel \}$. Also, due to (\ref{91}) and the dominated convergence theorem (see, e.g., p. 49 in \cite{Ash}),
\begin{equation}\lim_{N\rightarrow\infty}\int^t_0\parallel f(u^N(y^*(t^{\prime\ })),y^*(t^{\prime\ }))-f(u^*(t^{\prime\ }),y^*(t^{\prime\ }))\parallel dt^{\prime\ }=0. \label{96} \end{equation}
Let us introduce the notation
\begin{center}$ k_t(N)\stackrel{def}{=} L_2 meas\{A_t(N)\}+\int^t_0\parallel f(u^N(y^*(t^{\prime\ })),y^*(t^{\prime\ }))-f(u^*(t^{\prime\ }),y^*(t^{\prime\ }))\parallel dt^{\prime\ }$\end{center}
and rewrite the inequality (\ref{94}) in the form
\begin{equation} \parallel y^N(t)-y^*(t)\parallel\leq L_1\int^t_0\parallel y^N(t^{\prime\ })-y^*(t^{\prime\ })\parallel dt^{\prime\ } + k_t(N), \label{97} \end{equation}
which, by the Gronwall-Bellman lemma (see, e.g., p.218 in \cite{Bardi97}), implies that
\begin{equation} \max_{t\in [0,T^*]}\parallel y^N(t)-y^*(t)\parallel\leq k_t(N)e^{L_1T^*}.  \label{98} \end{equation}
Since, by (\ref{66}) and (\ref{96}),
\begin{equation} \lim_{N\rightarrow\infty}k_t(N)=0,     \label{99} \end{equation}
(\ref{98}) implies (\ref{68}).

For any $t\in [0,T^*]$ such that $u^N(\cdot)$ is Lipschitz continuous on $Q_t$, one has
\begin{center} $\parallel u^N(y^N(t))-u^*(t)\parallel \leq\parallel u^N(y^N(t))-u^N(y^*(t))\parallel+\parallel u^N(y^*(t))-u^*(t)\parallel$   \end{center}
\begin{center} $\leq L\parallel y^N(t)-y^*(t)\parallel+\parallel u^N(y^*(t))-u^*(t)\parallel.$ \end{center}
 The latter implies (\ref{67}) (due to (\ref{98}), (\ref{99}) and due to (\ref{91})).
 To finalize the proof, note that from (\ref{approx-u-per-extra}) it follows that
 $$
  \left|\frac{1}{T^*}\int^{T^*}_0 g(u^N(t),\tilde{y}^N(t))dt-G^*\right|
 $$
 $$
= \left|\frac{1}{T^*}\int^{T^*}_0 g(u^N(t),\tilde{y}^N(t))dt - \frac{1}{T^*}\int^{T^*}_0 g(u^*(t),y^*(t))dt\right|
 $$
$$
\leq \frac{1}{T^*}\int^{T^*}_0 ||g(u^N(t),\tilde{y}^N(t))- g(u^N(t),y^N(t))||dt
$$
$$
+ \frac{1}{T^*}\int^{T^*}_0 ||g(u^N(t),y^N(t))dt -  g(u^*(t),y^*(t))||dt
$$
$$
\leq \frac{L}{T^*}\int^{T^*}_0 [||\tilde{y}^N(t)- y^N(t)|| + || y^N(t) - y^*(t)|| + || u^N(t) - u^*(t)||] dt,
$$
where $L$ is a Lipschitz constant. The latter implies (\ref{102}) (due to (\ref{67}), (\ref{68}) and (\ref{approx-u-per-extra})). $ \ \Box$

 \textbf{Proof of Proposition \ref{Main-3}}. For any $t\in [0,T^*]$, one has
\begin{center}$\parallel \tilde{y}^N(t)-y^N(t)\parallel \leq \parallel \tilde{y}^N(0)-y^N(0)\parallel +\int^t_0 \parallel f(u^N(t^{\prime }),\tilde{y}^N(t^{\prime\ }))-f(u^N(t^{\prime\ }), y^N(t^{\prime\ }))\parallel\leq. $\end{center}
\begin{center}$ \leq \parallel \tilde{y}^N(0)-y^N(0)\parallel + L\int^t_0 \parallel \tilde{y}^N(t^{\prime\ })-y^N(t^{\prime\ })\parallel dt^{\prime\ }, $ \end{center}
which, by the Gronwall-Bellman Lemma, implies that
\begin{equation} \max_{t\in [0,T^*]}\parallel \tilde{y}^N(t)-y^N(t)\parallel \leq \parallel \tilde{y}^N(0)-y^N(0)\parallel e^{LT^*}. \label{103}\end{equation}
Due to  Assumption IV  and the periodicity condition $\tilde{y}^N(0)=\tilde{y}^N(T^*)$,  the following relationships are valid:
\begin{center}$ \parallel \tilde{y}^N(0)-y^N(0)\parallel \leq \parallel \tilde{y}^N(0)-y^N(T^*)\parallel +\parallel y^N(T^*)-y^N(0)\parallel  $ \end{center}
\begin{center}$ = \parallel \tilde{y}^N(T^*)-y^N(T^*)\parallel +\parallel y^N(T^*)-y^N(0)\parallel  $ \end{center}
\begin{center}$\leq \xi(T^*)\parallel\tilde{y}^N(0)-y^N(0)\parallel+\parallel y^N(T^*)-y^N(0)\parallel. $ \end{center}
Note that $y^N(0)=y^*(0)=y^*(T^*) $. Hence (see also (\ref{68})),
$$
\parallel y^N(T^*)-y^N(0)\parallel = \parallel y^N(T^*)-y^*(T^*)\parallel \leq \nu(N)
$$
\begin{center}$\Rightarrow \ \ \ \ \ \ \parallel \tilde{y}^N(0)-y^N(0) \parallel \leq \xi(T^*)\parallel \tilde{y}^N(0)-y^N(0)\parallel + \nu(N)$ \end{center}
\begin{center}$\Rightarrow \ \ \ \ \ \ \displaystyle \parallel \tilde{y}^N(0)-y^N(0)\parallel \leq \frac{\nu(N)}{1-\xi(T^*)}.$ \end{center}
Substituting the above inequality into (\ref{103}) one obtains
\begin{center}$\displaystyle \max_{t\in[0,T^*]}\parallel\tilde{y}^N(t)-y^N(t)\parallel\leq \frac{\nu(N)}{1-\xi(T^*)}e^{LT^*}.$ \end{center}
This proves (\ref{e-new-3}). The validity of (\ref{102}) is established as above. $\ \Box $

\section{Numerical example (swinging a nonlinear pendulum)}
\label{sec:6}
Consider the problem of periodic optimization of the nonlinear pendulum
\begin{equation}\label{e-pend-1}
x''(t)+0.3 x'(t) + 4 \sin(x(t)) = u(t)
\end{equation}
with the controls being restricted by the inequality $ \ |u(t)|\leq 1$ and with
the objective function being of the form
\begin{equation}\label{e-pend-2}
\inf_{u(\cdot), T}\frac{1}{T}\int_0^T (u^2(t)-x^2(t))dt.
\end{equation}
By re-denoting $x(t) $ and $x'(t)$ as $y_1(t)$ and $y_2(t) $ respectively, the above problem is reduced to  a special case of the periodic optimization problem (\ref{12}) with
$$
y=(y_1,y_2),  \ \ \ f(u,y)=(f_1(u,y),f_2(u,y))\stackrel{def}{=}(y_2,u-0.3y_2-4\sin(y_1)),  \ \ \  \  $$
$$g(u,y)\stackrel{def}{=}u^2-y_1^2 $$
and with
$$ U\stackrel{def}{=}[-1,1]\in R^1, \   \   \   \  \ Y\stackrel{def}{=}\{(y_1,y_2) \mid y_1 \in [-1.7, 1.7], \  \ y_2\in [-4,4]\} \in R^2 $$
(note that the set $Y$ is chosen to be large enough to contain all periodic solutions of the system under consideration).

The SILP problem (\ref{71}) was formulated for this problem with the use of the monomials $\phi_{i_1,i_2}(y)\stackrel{def}{=}y_1^{i_1}y_2^{i_2}, \  \  \ i_1,i_2=0,1,...,J$, as the functions $\phi_{i}(\cdot)$ defining $W_N(y_0)$ in (\ref{72}). Note that in this case the number $N$ in (\ref{72}) is equal to $(J+1)^2-1$. This problem and its dual were solved with the algorithm proposed in \cite{GRT} for the case $J=10 ( N= 120)$.
 In particular, the coefficients $\lambda^N_{i_1,i_2}$ defining the optimal solution of the corresponding  $N$-approximating maximin problem
\begin{equation} \psi^N(y)=\sum_{0<i_1+i_2\leq 10} \lambda ^N_{i_1,i_2}y_1^{i_1}y_2^{i_2} \end{equation}
were found (note the change of notations with respect to (\ref{49})), and the optimal value of the SILP was evaluated to be $\approx -1.174$.

In this case the problem (\ref{u-N}) takes  the form
$$
\min_{u\in [-1,1]}\{\frac{\partial\psi^N(y_1, y_2)}{\partial y_1}y_2 + \frac{\partial\psi^N(y_1, y_2)}{\partial y_2}(u-0.3y_2-4sin(y_1))+ (u^2-y_1^2) \}.
$$
The solution of the latter leads to the following representation for
 $u^N(y)$:
\begin{equation} u^N(y)=\left\{\begin{array}{ccc}
                            -\frac{1}{2}\frac{\partial\psi^N(y_1, y_2)}{\partial y_2} & {\rm if} & |\frac{1}{2}\frac{\partial\psi^N(y_1, y_2)}{\partial y_2}| \leq 1, \\
                            -1 & {\rm if} & -\frac{1}{2}\frac{\partial\psi^N(y_1, y_2)}{\partial y_2}<-1, \\
                            1 & {\rm if} & -\frac{1}{2}\frac{\partial\psi^N(y_1, y_2)}{\partial y_2} >1.
                          \end{array} \right.
 \end{equation}
 Substituting this control into the system (\ref{1}) and integrating it with the ode45 solver of MATLAB allows one to obtain the periodic ($T^*\approx 3.89$) state trajectory $\tilde{y}^N(t)=(\tilde{y}_1^N(t),\tilde{y}_2^N(t))$ (see Figure 1) and the control trajectory $u^N(t)$ (see Figure 2). The value of the objective function numerically evaluated on the state control trajectory thus obtained is $\approx -1.174$, the latter being the same as in SILP (within the given proximity).  Note that the  marked dots in Fig. 1 correspond to the concentration points of the measure $\gamma^N$ (see (\ref{78})) that solves (\ref{71}).
 The fact that the obtained state trajectory passes near these points and, most importantly, the fact that the value of the objective function  obtained via integration  is the same (within the given proximity) as the optimal value of the SILP problem indicate that the
  admissible solution found is a good approximation of the optimal one.

\begin{figure}
\begin{minipage}[b]{0.5\linewidth}
\begin{center}
\hspace*{-20mm}
\includegraphics[width=65mm]{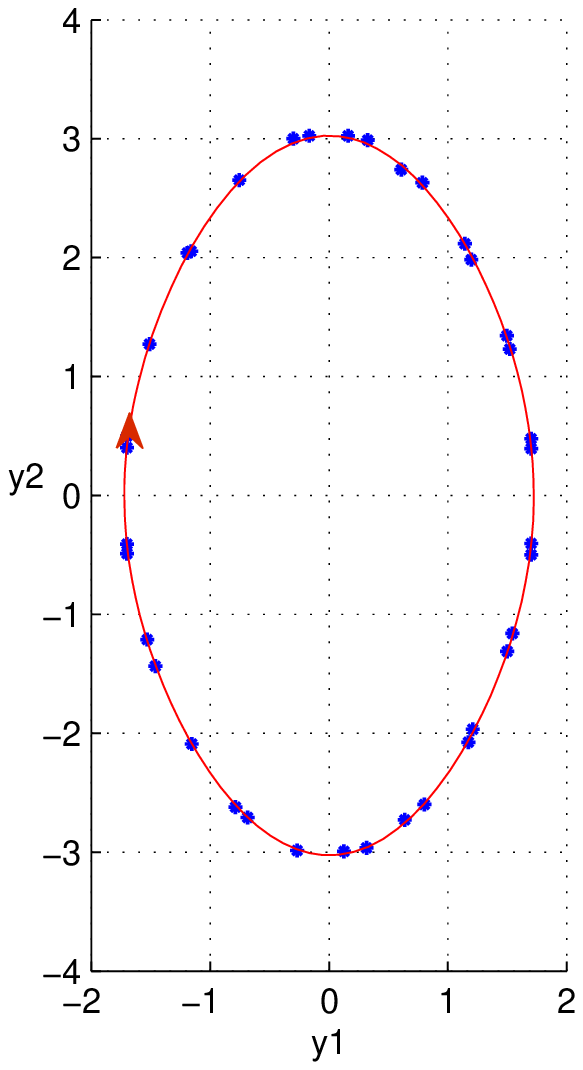} \\[3mm]
\hspace*{-1cm}{Fig.1: Near optimal state trajectory}
\end{center}
\end{minipage}
\begin{minipage}[b]{0.4\linewidth}
\begin{center}
\hspace*{-10mm}
\includegraphics[width=65mm]{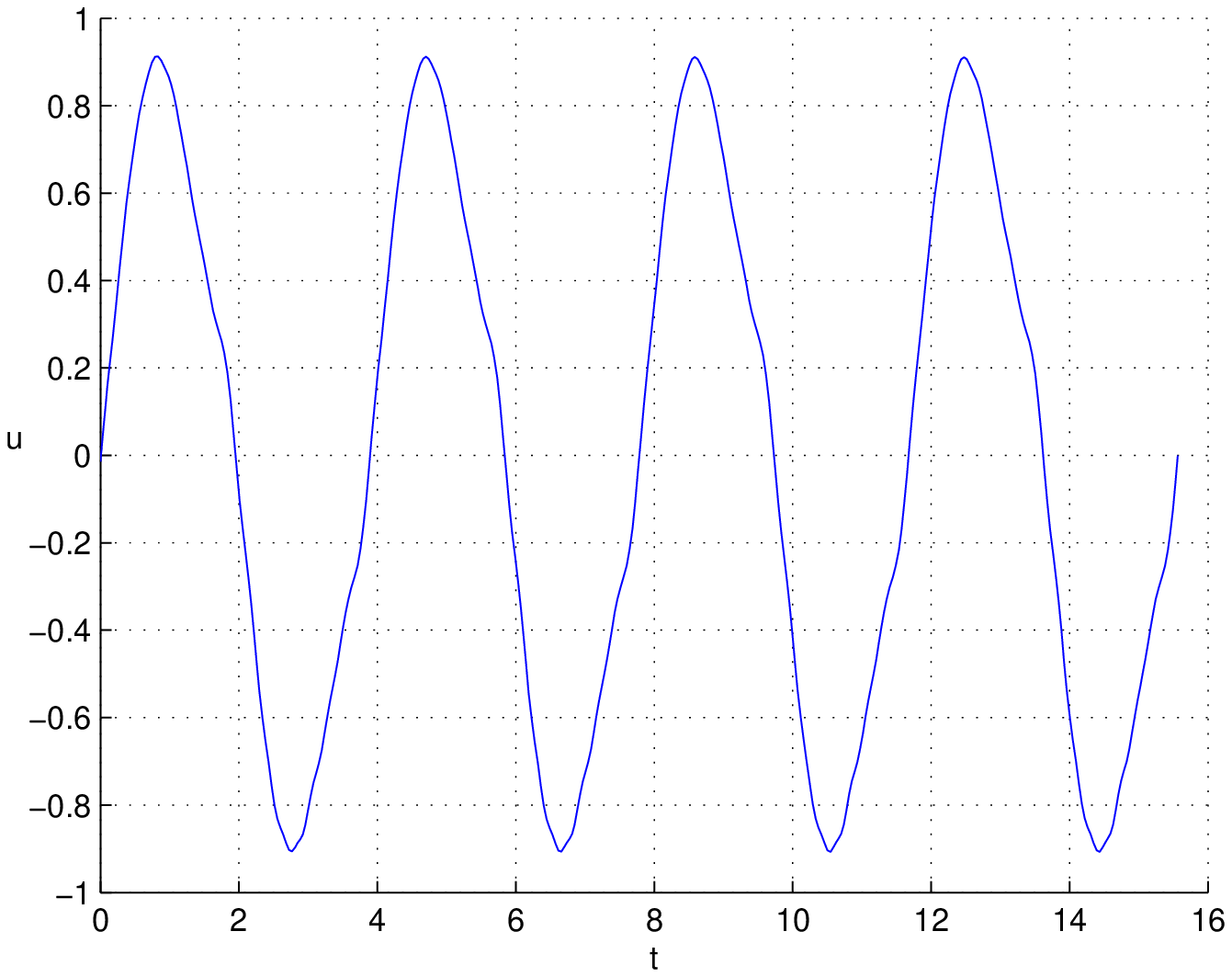} \\[3mm]
\hspace*{-1cm}{Fig.2: Near optimal control trajectory}
\end{center}
\end{minipage}
\end{figure}

\section{Algorithm for numerical solution of the SILP problem}

In this section, we describe an algorithm for solving the SILP problem (\ref{71}). It finds the optimal solution $\gamma ^N$ of the SILP problem  via solving a sequence of finite dimensional LP problems, each time augmenting the grid of Dirac measures concentration points with a new point found as an optimal solution of a certain nonlinear optimization problem. Note that, as it was shown in previous section, the solution of the latter  can be used for numerical construction of a near optimal control in the periodic optimization problem (\ref{12}).

For simplicity, let us denote $ X \stackrel{def}{=}U\times Y $, $x \stackrel{def}{=}(u,y)$ and $h_i(x)\stackrel{def}{=}\triangledown \phi_i(y)^Tf(u,y)$. Thus our SILP problem $G^*_N$  can be rewritten as follows
\begin{equation} \min_{\gamma\in W_N}\int_{U\times Y}g(x)\gamma (dx) \stackrel{def}{=}G^*_N, \label{71*} \end{equation}
where
 \begin{equation} W_N \stackrel{def}{=} \{  \gamma \in P(X): \ \int_{X}h_i(x)\gamma (dx)=0, \ \ i=1,...,N \}. \label{72.x} \end{equation}

Let points $ x_l \in X, \  \ l=1,...,M \  \ $ (note that $M \geq N+1$) be chosen to define an initial grid  $\Omega_0$ on $X$. That is
\begin{center}$ \Omega_0 = \{ x_l \in X, \  \  \ l=1,...,M \}.  $\end{center}
At every iteration a new point is defined and added to this set. Assume that after $k$ iterations, the points $x_{M+1},..., x_{M+k}$ have been defined and the set $\Omega_k$ has been constructed. Namely,
\begin{center} $ \Omega_k = \{ x_l \in X, \  \  \ l=1,..., M+k \}. $  \end{center}
The iteration $k+1 \ \ (k=0,1,...)$ is described as follows:

\begin{enumerate}[(i)]

\item Find an optimal solution of the LP problem
\begin{equation} \displaystyle \min_{\gamma \in  W_{\Omega_k}} \{ \sum^{M+k}_{l=1} g(x_l)\gamma_l \} \stackrel{def}{=}G^k,  \label{105} \end{equation}
where
\begin{equation} \displaystyle W_{\Omega_k} \stackrel{def}{=} \{ \gamma : \gamma = \{\gamma_l\} \geq 0, \  \  \sum^{M+k}_{l=1}\gamma_l=1,  \  \  \sum^{M+k}_{l=1} h_i(x_l) \gamma_l=0, \  \  \ i=1,...,N \}.  \end{equation}
Also find an optimal solution $\lambda^k = (\lambda^k_0, \  \lambda^k_i, \  i=1,...,N)$ of the problem dual with respect to (\ref{105}).  The latter being of the form
\begin{equation} \max \{ \lambda_0 : g(x_l)+ \sum^N_{i=1} h_i(x_l)  \lambda_i  \geq \lambda_0, \  \  \forall l=1,...,M+k \}. \label{110} \end{equation}

\item Find an optimal solution $x_{M+k+1}$ of the problem
\begin{equation}\min_{x \in X} \{g(x) + \sum^N_{i=1} h_i(x)  \lambda^k_i \} \stackrel{def}{=} a^k, \label{106} \end{equation}
where $\lambda^k = (\lambda^k_0,\  \lambda^k_i, \  i=1,...,N)$ is an optimal solution of the problem (\ref{110}).

\item Define the set $\Omega_{k+1}$ by the equation
\begin{center} $ \Omega_{k+1}=\Omega_k \bigcup \  \{ x_{M+k+1} \}. $\end{center}

\end{enumerate}
Here and in what follows, $k$ stands for the number of an iteration. Note that, by construction
\begin{equation}\label{e-momotonicity} G^{k+1}\leq G^k, \  \  k=1,2,... \ . \end{equation}
In the next section, we establish that, under certain regularity (non-degeneracy) conditions, the optimal value $G^k$ of the problem (\ref{105})
converges (as $k$ tends to infinity) to the optimal value $G^*_N $ of the problem (\ref{71}) (see Theorem \ref{8.7*} below).

\section{Convergence of the algorithm}
\label{section continuity}

Let  $ \mathfrak{X}= \{ x_1, x_2, ..., x_{N+1}  \} \in X^{N+1} $ and let
$\gamma(\mathfrak{X})=\{ \gamma_j(\mathfrak{X}),  \ \forall j=1,...,N+1\}\geq 0  \label{8.3}  $ satisfy the system of $N+1 $ equations
\begin{equation} \sum^{N+1}_{j=1}\gamma_j(\mathfrak{X}) = 1, \    \    \    \sum^{N+1}_{j=1} h_i (x_j)  \gamma_j(\mathfrak{X})=0, \  \   \ \forall i=1,...,N. \label{8.1} \end{equation}
Assume that the solution of the system (\ref{8.1}) is unique (that is, the system is non-singular) and
define
\begin{equation}  G(\mathfrak{X}) \stackrel{def}{=} \sum^{N+1}_{j=1} g(x_j)\gamma_j(\mathfrak{X}).  \label{8.2} \end{equation}
Also, let $\lambda(\mathfrak{X}) = \{ \lambda_0(\mathfrak{X}), \lambda_1(\mathfrak{X}), ..., \lambda_N(\mathfrak{X}) \} $ be a solution of the system
\begin{equation} \lambda_0(\mathfrak{X}) - \sum^N_{i=1} h_i(x_j) \lambda_i(\mathfrak{X})  =g(x_j),  \ \ j=1,...,N+1. \label{8.4} \end{equation}
and let
\begin{equation} a(\mathfrak{X})\stackrel{def}{=} \min_{x\in X}\{g(x)+\sum_{i=1}^N h_i(x) \lambda_i(\mathfrak{X})  \}.  \label{8.5} \end{equation}

\begin{lemma}\label{8.1*}
For any $\mathfrak{X} \subset X^{N+1}$,
\begin{equation} -\lambda_0(\mathfrak{X})+a(\mathfrak{X})\leq 0 . \label{8.7} \end{equation}
\end{lemma}

\textbf{Proof} By the definition of $\mu^*_N$ (see (\ref{45})), \  \  $ a(\mathfrak{X})\leq \mu^*_N $. \ \ Also, due to the duality theorem (see (\ref{79}))
\begin{equation}\mu^*_N = G^*_N. \label{8.8} \end{equation}
Hence,
\begin{equation} a(\mathfrak{X})\leq G^*_N.  \label{8.9} \end{equation}
Note that, by multiplying the $j^{th}$ equation in (\ref{8.4}) by $\gamma_j(\mathfrak{X})$ and by summing up the resulted equations over $j=1,...,N+1$ one can obtain (using (\ref{8.1}) and (\ref{8.2})) that
\begin{equation} \lambda_0(\mathfrak{X}) =\sum^{N+1}_{j=1} g(x_j)\gamma_j(\mathfrak{X}) =G(\mathfrak{X}). \label{8.10} \end{equation}
Since $G(\mathfrak{X})\geq G^*_N$, from (\ref{8.10}) it follows that $ \lambda_0(\mathfrak{X})\geq G^*_N$. The later and (\ref{8.9}) proves (\ref{8.7}).  $\Box$

\begin{corollary}
\label{8.2*}
If
\begin{equation} -\lambda_0(\mathfrak{X})+a(\mathfrak{X})= 0, \label{8.11} \end{equation}
then $\gamma(\mathfrak{X}) \stackrel{def}{=}\sum^{N+1}_{j=1}\gamma_j(\mathfrak{X}) \delta_{x_j}$ (where $\delta_{x_j} $ is
the Dirac measure concentrated at $x_j $) is an optimal solution of the SILP problem (\ref{71})
and $\psi(y)\stackrel{def}{=}\sum^N_{i=1}\lambda_i(\mathfrak{X})\phi_i(y)$ is an optimal solution of the N-approximating problem (\ref{45}).
\end{corollary}

\textbf{Proof} Due to (\ref{8.9}) and (\ref{8.10}) we have
\begin{equation} a(\mathfrak{X})\leq \mu^*_N =G^*_N \leq G(\mathfrak{X}) =\lambda_0(\mathfrak{X}).   \label{8.12} \end{equation}
From (\ref{8.11}) and (\ref{8.12}) it follows that
\begin{equation} a(\mathfrak{X})= \mu^*_N \label{8.13} \end{equation}
and \begin{equation} G(\mathfrak{X})=G^*_N. \label{8.14}\end{equation}
The equalities (\ref{8.13}) and (\ref{8.14}) proves the statements of the corollary \ref{8.2*}. $\Box$

\begin{definition} \label{8.3*} An $(N+1)$-tuple $\mathfrak{X} \subset X^{N+1}$ is called \emph{regular} if
\begin{enumerate} [(i)]
\item the system (\ref{8.1}) has a unique  solution (that is, the corresponding $(N+1)\times (N+1)$ matrix is not singular  ); and

\item the solution $\gamma(\mathfrak{X})$ of this system is positive. That is,
\begin{center}$\gamma(\mathfrak{X})\stackrel{def}{=}\{\gamma_j(\mathfrak{X}), \ \ \forall j=1,...,N+1 \} > 0 .$ \end{center}
\end{enumerate}
\end{definition}

Assume that $\mathfrak{X}$ is regular and  examine the case when
\begin{equation} -\lambda_0(\mathfrak{X}) + a(\mathfrak{X})< 0. \label{8.15} \end{equation}
Denote by $A(\mathfrak{X})$ the $(N+1)\times (N+1)$ matrix of the system (\ref{8.1}). That is
\begin{center} $A(\mathfrak{X})=\{ H(x_1), H(x_2),..., H(x_{N+1})  \},$ \end{center}
where columns $H(x_j), j=1,2,...,N+1$ are defined as follows
\begin{center} $ H(x_j)=(1,h_1(x_j), h_2(x_j), ..., h_N(x_j))^T $. \end{center}
Using this notations, the solution $\gamma(\mathfrak{X})$ of the system of equations (\ref{8.1}) will be written in the form
\begin{equation} \gamma(\mathfrak{X})=A^{-1}(\mathfrak{X})b, \label{8.16} \end{equation}
where $b=(1,0,0,...,0)^T$.
Similarly, the solution of the system (\ref{8.4}) is defined as follows
\begin{equation} \lambda(\mathfrak{X})\stackrel{def}{=}(\lambda_0(\mathfrak{X}),\  \lambda_i(\mathfrak{X}), \  \   i=1,...,N)= (A^T(\mathfrak{X}))^{-1} c(\mathfrak{X}), \label{8.17} \end{equation}
where $c(\mathfrak{X})\stackrel{def}{=} (\ g(x_1),\ g(x_2),\ ...\ ,\ g(x_{N+1}) \ )^T$.

Let  $\mathcal{B}(\mathfrak{X})$ stand for the set of the optimal solutions of the problem (\ref{8.5}). That is,
\begin{equation} \mathcal{B}(\mathfrak{X})\stackrel{def}{=} Arg \min_{x\in X} \{ g(x) + \sum^N_{i=1}h_i(x) \lambda_i(\mathfrak{X})  \}. \label{8.6} \end{equation}
Choose an arbitrary $\tilde{x} \in  \mathcal{B}(\mathfrak{X})$ and consider the system of equations
 \begin{equation} \sum^{N+1}_{j=1}\gamma_j(\theta)+\theta=1, \  \  \  \sum^{N+1}_{j=1}h_i(x_j)\gamma_j(\theta) + h_i(\tilde{x})\theta=0, \  \  \forall i=1,...,N ,  \label{8.18} \end{equation}
where $\theta \geq 0$ is a parameter and $\gamma_j(\theta), j=1,2,...,N+1$ are defined by the value of this parameter. Note that this system also can be presented in the form
\begin{equation}  A(\mathfrak{X})\gamma(\theta)+ H(\tilde{x})\theta =b. \label{e-extra}\end{equation}
Let also $G(\mathfrak{X},\theta)$ be defined by the equation
\begin{equation} G(\mathfrak{X}, \theta) \stackrel{def}{=}\sum^{N+1}_{j=1}g(x_j)\gamma_j(\theta)+g(\tilde{x})\theta  .  \label{8.19} \end{equation}
Observe that, by multiplying the $j^{th}$ equation in (\ref{8.4}) by $\gamma_j(\theta)$ and summing up the  resulted equations over $j=1,...,N+1$ one can obtain (using (\ref{8.18}) and (\ref{8.19}))
\begin{equation} \lambda_0(\mathfrak{X})(1-\theta) +\sum_{i=1}^N h_i(\tilde{x})\lambda_i(\mathfrak{X}) \theta = G(\mathfrak{X}, \theta)-g(\tilde{x})\theta. \label{8.20} \end{equation}
Due to  (\ref{8.5}), (\ref{8.10}) and the fact that $\tilde{x}\in  \mathcal{B}(\mathfrak{X})$, from(\ref{8.20}) it follows that
\begin{equation} G(\mathfrak{X}, \theta) = G(\mathfrak{X}) +(-\lambda_0(\mathfrak{X}) + a(\mathfrak{X}))\theta.   \label{8.21} \end{equation}
By (\ref{8.15}), the bigger value of $\theta$ is, the better is the resulted value of the objective
function $ G(\mathfrak{X}, \theta)$. The value of $\theta $ is, however,  bounded from above  by  $\displaystyle \theta(\mathfrak{X},\tilde{x})\stackrel{def}{=}\min_{j,\ d_j>0} \frac{\gamma_j(\mathfrak{X})}{d_j(\mathfrak{X},\tilde{x})}$, where $\ d(\mathfrak{X},\tilde{x})\stackrel{def}{=}A^{-1}(\mathfrak{X})H(\tilde{x})$ (this constraint being implied by (\ref{e-extra})).
  Thus,  the improvement, induced by replacing one of the columns of $A(\mathfrak{X}) $  with the  column $H(\tilde{x})$ in accordance with Simplex Method (see \cite{Dantzig})  is determined by the following expression
\begin{equation}  [ -\lambda_0(\mathfrak{X}) +a(\mathfrak{X}) ] \ \theta(\mathfrak{X},\tilde{x}). \label{8.22} \end{equation}
Let us denote $\displaystyle \min_{j} \gamma_j(\mathfrak{X})\stackrel{def}{=} \eta (\mathfrak{X})>0$ and $\ \ \displaystyle \max_{j}d_j(\mathfrak{X},\tilde{x}) \stackrel{def}{=}\beta(\mathfrak{X},\tilde{x})>0$.
Also, let \begin{center} $\displaystyle \max_{\tilde{x}\in \mathcal{B}(\mathfrak{X})}\beta(\mathfrak{X},\tilde{x})\stackrel{def}{=}\beta(\mathfrak{X})$.\end{center}
Define
\begin{equation}\label{e-extra-2}\mathcal{V}(\mathfrak{X})\stackrel{def}{=}(-\lambda_0(\mathfrak{X})+a(\mathfrak{X}))\frac{\eta(\mathfrak{X})}{\beta(\mathfrak{X})}. \end{equation}

\begin{lemma} \label{8.4*}
If $\mathfrak{X}$ is regular and $a(\mathfrak{X})< \lambda_0(\mathfrak{X})$,  then $\forall \tilde{x}\in  \mathcal{B}(\mathfrak{X})$ the replacement of the one of the column in $A(\mathfrak{X})$ by $H(\tilde{x})$ (according to  Simplex Method) leads to the improvement of the objective value no less then $\mathcal{V}(\mathfrak{X})$. That is,
\begin{equation} (-\lambda_0(\mathfrak{X})+a(\mathfrak{X}))\theta(\mathfrak{X},\tilde{x})\geq \mathcal{V}(\mathfrak{X}). \label{8.23} \end{equation}
\end{lemma}

\textbf{Proof}. The proof is obvious. $\Box$

Note that if $\mathfrak{X}$ is regular, then any $\mathfrak{X}'\in \mathcal{D}_r\BYDEF \{\mathfrak{X}'\ : \ \parallel\mathfrak{X}'-\mathfrak{X}\parallel < r\} $ (where  $r>0$ is small enough) will be regular as well, with $\eta(\mathfrak{X}')$ being continuous function of $\mathfrak{X}'$  and $\beta(\mathfrak{X}',\tilde{x})$ being  continuous function of $\mathfrak{X}'$ and $\tilde{x}$.

\begin{lemma} \label{8.5*} Let $\mathfrak{X}$ be  regular. Then the function  $\beta(\cdot)$ is upper semicontinuous at $\mathfrak{X}$. That is,
\begin{equation} \varlimsup_{\mathfrak{X}^{l}\rightarrow \mathfrak{X}}\beta(\mathfrak{X}^{l})\leq \beta(\mathfrak{X}). \label{8.24} \end{equation}
\end{lemma}

\textbf{Proof.} 
Let us first of all show that
\begin{equation} \displaystyle \varlimsup_{\mathfrak{X}^l\rightarrow \mathfrak{X}} \mathcal{B}(\mathfrak{X}^l) \subset  \mathcal{B}(\mathfrak{X}). \label{8.25} \end{equation}
That is,  if $\forall x^l \in  \mathcal{B}(\mathfrak{X}^l)$ and  $\displaystyle \lim_{l\rightarrow \infty} x^l=x\in X$, then
\begin{equation} x \in  \mathcal{B}(\mathfrak{X}). \label{8.26} \end{equation}
The fact that $x^{l}\in  \mathcal{B} (\mathfrak{X}^l)$ means that
\begin{equation} g(x^l)+\sum^{N}_{i=1}\lambda_i(\mathfrak{X}^l)h_i(x^l) =a(\mathfrak{X}^l). \label{8.27}\end{equation}
Passing to the limit  as $l\rightarrow \infty$ in  (\ref{8.27}) one can obtain (having in mind that $\lambda(\cdot)$ and $a(\cdot)$ are continuous function in a neighbourhood of $\mathfrak{X}$), one obtains the equality
\begin{equation} g(x)+\sum^{N}_{i=1} \lambda_i(\mathfrak{X}) h_i(x)=a(\mathfrak{X}). \label{8.28}\end{equation}
The latter proves (\ref{8.26}) and thus establishes validity of (\ref{8.25}).
To prove (\ref{8.24}) let $\mathfrak{X}^l \rightarrow \mathfrak{X}$ as $l\rightarrow \infty$. \ \ Recall that
\begin{equation} \beta(\mathfrak{X}^l)=\max_{\tilde x\in  \mathcal{B}(\mathfrak{X}^l)} \beta(\mathfrak{X}^l,\tilde x) . \label{8.29} \end{equation}
Let $\tilde{x}^l \in \mathcal{B}(\mathfrak{X}^l)$ be such that maximum in (\ref{8.29}) is reached. Without loss of generality, one may assume that there exists
a limit
\begin{center} $ \displaystyle \lim_{l\rightarrow \infty} \tilde{x}^l = \tilde{x} \in \mathcal{B} , $ \end{center}
with the incusion being due to (\ref{8.25}). Thus, passing to the limit as $l\rightarrow \infty $ in (\ref{8.29}) one can get
\begin{equation} \displaystyle \lim_{l\rightarrow \infty } \beta(\mathfrak{X}^l) = \lim_{l\rightarrow \infty } \beta(\mathfrak{X}^l, \tilde{x}^l) =\beta(\mathfrak{X}, \tilde{x}) \leq \beta(\mathfrak{X}).  \end{equation}
\ $\Box$

\begin{lemma}\label{8.6*}
If $\mathfrak{X}$ is regular and $\ a(\mathfrak{X}) < \lambda_0(\mathfrak{X})\ $ then $\ \exists r_0 >0$ such that
\begin{equation} \mathcal{V}(\mathfrak{X}^{'})\geq \frac{\mathcal{V}(\mathfrak{X})}{2}, \label{8.30} \end{equation}
for any $\mathfrak{X}'$  such that $\   ||\mathfrak{X}^{'}-\mathfrak{X}||\leq r_0 $,
where $\mathcal{V}(\cdot) $ is defined by (\ref{e-extra-2}).
\end{lemma}

\textbf{Proof}  To prove (\ref{8.30}) let us show that $\mathcal{V}(\cdot)$ is  lower semicontinuous at $\mathfrak{X}$.
According to (\ref{e-extra-2}),
\begin{equation} \mathcal{V}(\mathfrak{X}^{'})\stackrel{def}{=}(-\lambda_0(\mathfrak{X}^{'})+a(\mathfrak{X}^{'})) \frac{ \eta (\mathfrak{X}^{'}) }{ \beta(\mathfrak{X}^{'}) }.  \label{8.31} \end{equation}
By taking the lower limit in (\ref{8.31}) one can obtain (using Lemma \ref{8.5*})
\begin{center}$ \displaystyle \varliminf_{\mathfrak{X}^{'}\rightarrow \mathfrak{X}}\mathcal{V}(\mathfrak{X}^{'})= (-\lambda_0(\mathfrak{X})+a(\mathfrak{X})) \  \eta(\mathfrak{X}) \varliminf_{\mathfrak{X}^{'}\rightarrow \mathfrak{X}} \frac{1}{\beta(\mathfrak{X}^{'})}$\end{center} \begin{center}$ = (-\lambda(\mathfrak{X})+a(\mathfrak{X})) \  \eta(\mathfrak{X}) \frac{1}{\displaystyle \varlimsup_{\mathfrak{X}^{'}\rightarrow \mathfrak{X}}\beta(\mathfrak{X}^{'})}.  $\end{center}
Due to (\ref{8.24}), it follows that
\begin{equation}  \varliminf_{\mathfrak{X}^{'}\rightarrow \mathfrak{X}}\mathcal{V}(\mathfrak{X}^{'})\geq (-\lambda(\mathfrak{X})+a(\mathfrak{X}))\frac{\eta(\mathfrak{X})}{\beta(\mathfrak{X})}= \mathcal{V}(\mathfrak{X}). \label{8.33}\end{equation}
Thus,  $\mathcal{V}(\mathfrak{X})$ is lower semicontinuous at $\mathfrak{X}$. Hence, for any $\epsilon >0$, there exists $r >0$ such that for
\begin{equation} \mathcal{V}(\mathfrak{X}^{'})\geq \mathcal{V}(\mathfrak{X})-\epsilon \label{8.34} \end{equation}
for $\mathfrak{X}^{'} $ such that $||\mathfrak{X}^{'}-\mathfrak{X}||\leq r$.
By taking  $\epsilon = \frac{\mathcal{V}(\mathfrak{X})}{2}$ in (\ref{8.34}), one establishes  the validity of (\ref{8.30}).
\ $\Box$

Let $\ \mathfrak{X}^k = \{x^k_1,..., x^k_{N+1}\}\in X^{N+1}, \ \ \ x^k_j\in \Omega_k, \ j=1,...N+1, \ \ $ and
$\gamma^k_j\geq 0, \ \ j=1,..., N+1, $  be such that
\begin{equation}\label{e-extra-3-1} \sum_{l=1}^{N+1}\gamma^k_jg(x^k_j)= G^k \ \end{equation}
and
\begin{equation}\label{e-extra-4-1}
 \sum_{j=1}^{N+1}\gamma^k_j =1, \ \ \ \  \ \ \  \sum_{j=1}^{N+1}h_i(x^k_j)\gamma^k_j=0, \ \ i=1,..., N.
 \end{equation}
That is, $\{\gamma^k_1,..., \gamma^k_{N+1}\}$
are basic components of an optimal basic solution of the problem (\ref{105}).

Let $\Lambda\subset X^{N+1}$ stand for the set of cluster (limit) points
of $\{ \mathfrak{X}^k \}, \ \ k=1,2,...$.

\begin{theorem} \label{8.7*}  Let there exists at least one  regular $\mathfrak{X}\in \Lambda $.
Then
\begin{equation}\label{e-Th-1}  \displaystyle \lim_{k\rightarrow \infty} G^k =G^*_N .  \end{equation}
Also, if $ \ \mathfrak{X}\stackrel{def}{=}\{ x_j \} \in \Lambda\ $ is regular, then $\ \sum^{N+1}_{j=1}\gamma_j(\mathfrak{X})\delta_{x_j}\ $ is an optimal solution of the SILP problem (\ref{71}) and  $\sum^N_{i=1}\lambda_i(\mathfrak{X})\phi_i(y)$ is an optimal solution of the N-approximating problem (\ref{45}), where $\ \gamma(\mathfrak{X})=\{ \gamma_j(\mathfrak{X}),  \ \forall j=1,...,N+1\}$
is the solution of the system (\ref{8.1}) and $ \ \lambda(\mathfrak{X}) =\{\lambda_0(\mathfrak{X}), \lambda_1(\mathfrak{X}),...,\lambda_{N+1}(\mathfrak{X}) \} $ is the solution of the system (\ref{8.4}).

\end{theorem}

\textbf{Proof.} \ \ Let $\mathfrak{X}\in \Lambda $ be regular. By definition, there exists a subsequence $\{ \mathfrak{X}^{k_l} \} \in \{ \mathfrak{X}^k \}$ such that \begin{equation}\displaystyle \lim_{l\rightarrow \infty}  \mathfrak{X}^{k_l}=\mathfrak{X}.\label{8.33*} \end{equation} Thus,
\begin{equation} \lim_{l\rightarrow \infty} G(\mathfrak{X}^{k_l})=G(\mathfrak{X}).  \label{8.35} \end{equation}
Also,  due to (\ref{e-momotonicity}), there exists a limit
\begin{equation} \lim_{k\rightarrow \infty}G^k \stackrel{def}{=}\widetilde{G}. \label{8.36}\end{equation}
Since, by definition, $G(\mathfrak{X}^{k_l}) = G^{k_l} $, it follows that
\begin{equation} G(\mathfrak{X})=\widetilde{G}. \label{8.37-1}\end{equation}
By Lemma 8.1 there are two possibilities
\begin{enumerate}[(i)]
\item  $-\lambda_0(\mathfrak{X})+a(\mathfrak{X})=0$, and

\item  $-\lambda_0(\mathfrak{X})+a(\mathfrak{X})<0$.
\end{enumerate}

If $a(\mathfrak{X})=\lambda_0(\mathfrak{X})$, then  by (\ref{8.14}),  $G(\mathfrak{X})=G_N^*$ and, hence, the validity
of (\ref{e-Th-1}) follows from
(\ref{8.36}) and (\ref{8.37-1}).
Let us prove that (ii) leads to a contradiction. Observe that, due to (\ref{8.36}),
\begin{equation} |G^{k_l+1}-G^{k_l}|< \frac{\mathcal{V}(\mathfrak{X})}{4} \label{8.38} \end{equation}
for $l$ large enough. On the other hand, due to (\ref{8.33*})
 \begin{center}$\|\mathfrak{X}^{k_l}-\mathfrak{X}\|< r_0, $ \end{center}
for $l$ large enough, with $r_0$ being as in Lemma \ref{8.6*}. Hence,  by this lemma and by (\ref{e-momotonicity}),
\begin{equation} |G^{k_l+1}-G^{k_l}|\geq \frac{\mathcal{V}(\mathfrak{X})}{2}. \label{8.39} \end{equation}
The latter contradicts (\ref{8.38}). Thus, the inequality $-\lambda_0(\mathfrak{X})+a(\mathfrak{X})<0$ can not be valid and
(\ref{e-Th-1}) is proved. Also, by Corollary \ref{8.2*}, from the fact that $-\lambda_0(\mathfrak{X})+a(\mathfrak{X})=0$
it follows that $\ \sum^{N+1}_{j=1}\gamma_j(\mathfrak{X})\delta_{x_j}\ $ is an optimal solution of the SILP problem (\ref{71}) and  $\sum^N_{i=1}\lambda_i(\mathfrak{X})\phi_i(y)$ is an optimal solution of the N-approximating problem (\ref{45}).

\ $\Box$

REMARK. Note that the assumption that there exists at least one  regular $\mathfrak{X}\in \Lambda$ is much weaker
than one used in proving the convergence of a similar algorithm in \cite{GRT}, where it was assumed  
that the optimal solutions of the LP problems (\ref{105}) are \lq\lq uniformly non-degenerate" (that is, they remain greater than some given positive number for $k=1,2,...$; see Proposition 6.2 in \cite{GRT}).

\vspace{-0.5cm}
\footnotesize{

\end{document}